\documentclass[11pt,a4paper]{article}

\usepackage{graphicx}
\usepackage{amsmath,amsthm}
\usepackage[numbers]{natbib}
\usepackage{microtype}
\usepackage[bitstream-charter]{mathdesign}
\usepackage[final]{showkeys}

\numberwithin{equation}{section}
\allowdisplaybreaks

\newtheorem{Theorem}{Theorem}[section]

\title{A note on fractional linear pure birth and pure death processes in epidemic models}

\author{$\text{Roberto Garra}_1$, $\text{Federico Polito}_2$\\
	\footnotesize (1) -- Dipartimento di Fisica, ``Sapienza'' Universit\`a di Roma\\
	\footnotesize Piazzale Aldo Moro 5, 00185 Rome, Italy.\\
	\footnotesize Email address: rolinipame@yahoo.it \\
	\footnotesize (2) -- Dipartimento di Matematica, Universit\`a degli studi di Roma
	``Tor Vergata''\\
	\footnotesize Via della Ricerca Scientifica 1, 00133, Rome, Italy.\\
	\footnotesize Tel: +39-06-2020568, fax: +39-06-20434631\\
	\footnotesize Email address: polito@nestor.uniroma2.it
	}

\begin{document}

	\maketitle
	
	\begin{abstract}
		\noindent In this note we highlight the role of fractional linear birth and linear
		death processes recently studied in
		\citet{sakhno} and \citet{pol}, in relation to epidemic models with empirical power law
		distribution of the events.  Taking inspiration from a formal analogy  between the equation of self
		consistency of the epidemic type aftershock sequences (ETAS) model, and the fractional
		differential equation describing the mean value of fractional linear growth processes, we show some
		interesting applications of fractional modelling to study \textit{ab initio} epidemic processes
		without the assumption of any empirical distribution. We also show that, in the frame of fractional
		modelling, subcritical regimes can be linked to linear fractional death processes and supercritical
		regimes to linear fractional birth processes.
		
		Moreover we discuss a simple toy model to underline the possible application of these
		stochastic growth models to more general epidemic phenomena such as tumoral growth.
		
		\smallskip
		
		\noindent Keywords: \emph{ETAS model, fractional branching, birth process, death process,
			Mittag--Leffler functions, Wiener--Hopf integral.}
	\end{abstract}

	\section{Introduction}
		In recent years there has been a growing interest in fractional calculus
		modelling in different fields of applied sciences, from rheology to biology
		(see for example \citet{podlubny,debnath,mainardi,diethelm}).
		It is well known that fractional derivatives are a good instrument to handle memory mechanisms,
		which are
		very useful in studying anomalous diffusion processes.
		Indeed, the fractional derivative in the sense of Caputo is
		defined by a convolution of a power law with the ordinary
		derivative of the function, thus underlining its importance as an
		instrument to describe processes with emerging power law
		distributions. Recent papers \citep{sakhno,pol} describe fractional
		processes of death and birth in which the integer-order derivative
		in the difference-differential equations governing the state
		probabilities of the related classical models, is replaced with the Caputo fractional derivative.
		In this note we apply this
		generalisation to the study of epidemic processes. In particular, we take inspiration from
		an interesting analogy between the self consistency equation 
		of the ETAS (epidemic type aftershock sequences) model,
		used with success in statistical seismology, and the equation governing the mean behaviour
		of the fractional linear death model.
		The ETAS model, based on epidemic branching processes,
		\citep{saichev}, is used to describe seismic sequences of shocks in the context of a stochastic model
		in which each aftershock
		may trigger other aftershocks. This model is based on empirical power law distributions emerging
		from the analysis of
		seismic catalogues, i.e., the modified Omori law on the time distribution of aftershocks and the
		Gutenberg--Richter law for the distribution of their magnitudes.
		These types of models are greatly studied in literature
		(see for example \citet{sornette,helm}) and used in assessing seismic risk \citep{helm2}.
		The self consistency equation of the ETAS model \citep{sornette} is the classical
		Wiener--Hopf integral equation, but with simple
		handling it can be considered to be a fractional integral equation similar to
		the equation governing the mean behaviour
		of the fractional linear death model.
		Using this analogy we suggest an interpretation of the fractional death and birth processes 
		in the framework of epidemic models with power law distributions (e.g.\ ETAS).
		The organization of this paper is the following:
		in Section \ref{sec2} we recall the fundamental results on fractional 
		death and fractional birth models; in Section \ref{analysis}
		we describe the analytic formulation of the ETAS model from which we took inspiration 
		for the discussion which follows in Section \ref{sec4}, that is the relations between the ETAS and
		fractional processes; in Section \ref{sec5}, in order to highlight the utility of this
		stochastic framework, we give an example of application of this
		formalism to describe tumoral growth;
		finally, in Section \ref{sec6}
		we discuss the main results obtained.

	\section{Fractional linear pure death and pure birth processes}
		
		\label{sec2}
		Recently, in \citet{sakhno}, section 2, the fractional linear death process was introduced and studied.
		The fractional linear death process is a generalisation of the well-known classical
		linear death process. Fractionality is obtained by
		replacing the integer-order derivative in the difference-differential equations governing the state
		probabilities of the classical model with the Caputo fractional derivative \eqref{caputo} of order $\nu
		\in (0,1)$. We exclude the case $\nu=1$ as it is trivial and coincides with the classical case.
		Let $n_0$ be the size of the population at time $0$ and let $M^\nu(t)$, $t>0$, $\nu \in (0,1)$,
		be the number of components in the population in the fractional linear
		death process at time $t$.
		
		The state probabilities $\text{Pr} \{ M^\nu(t) = k | M^\nu(0) = n_0 \}$,
		$0 \leq k \leq n_0$, solve the following Cauchy problem:
		\begin{align}
			\begin{cases}
				D^\nu p_k(t)=\mu (k+1) p_{k+1}(t)-\mu k p_{k}(t), & 0 \leq k \leq n_0, \: \nu \in (0,1), \\
				p_k ^{\nu}(0)=
				\begin{cases}
					1 & \text{if $k = n_0$}, \\
					0 & \text{if $k < n_0$},
				\end{cases}
			\end{cases}
		\end{align}
		where $\mu>0$ is the death intensity and in which we consider that $p_{n_0+1} (t)= 0$.
		
		We recall the following fundamental theorem \citep[page 73]{sakhno}:
		\begin{Theorem}
			The distribution of the fractional linear death process $M^{\nu}(t)$, $t>0$ with $n_0$ initial
			individuals and death rates $\mu_k=\mu k$, is given by:
			\begin{align}
				p^\nu_k(t) & = \text{Pr}\{M^\nu(t)=k|M^\nu(0)=n_0\} \\
				& = \binom{n_0}{k}\sum_{r=0}^{n_0-k} \binom{n_0-k}{r}(-1)^{r}
				E_{\nu,1}(-(k+r) \mu t^{\nu}) \notag
			\end{align}
			where  $0\leq k \leq n_0$, $t>0$, $\nu \in(0,1)$. The function $E_{\nu,1}(x)$ is the Mittag--Leffler
			function defined as
			\begin{align}
				E_{\nu,1}(x)=\sum_{k=0}^{+\infty}\frac{x^k}{\Gamma(\nu k+1)} \qquad x \in \mathbb{R}, \: \nu > 0.
			\end{align}
		\end{Theorem}
		In the same paper it is also shown that the mean value $\mathbb{E} M^\nu(t)$, $t>0$,
		of the fractional linear pure death process satisfies the fractional differential equation
		\begin{align}
			\label{ponte}
			\begin{cases}
				D^{\nu}\mathbb{E}M^{\nu}(t)=-\mu\mathbb{E}M^{\nu}(t), & t>0, \: \nu \in (0,1),\\
				\mathbb{E}M^{\nu}(0)=n_0,
			\end{cases}
		\end{align}
		which can be easily solved by means of the Laplace transform, obtaining 
		\begin{align}
			\mathbb{E}M^{\nu}(t)=n_0 E_{\nu,1}(-\mu t^{\nu}), \qquad t>0, \: \nu \in (0,1).
		\end{align}
		
		In \citet{pol}, the fractional pure birth process $\mathfrak{N}^\nu(t)$, $t>0$, $\nu \in (0,1)$,
		was studied with similar methods to those used for the fractional death.
		For the sake of our discussion, we recall only the fractional differential equation
		which governs the mean value $\mathbb{E}\mathfrak{N}^\nu(t)$, $t>0$, referring to the original paper
		for further details:
		\begin{align}
			\label{birth}
			\begin{cases}
				D^{\nu}\mathbb{E}\mathfrak{N}^{\nu}(t)
				= \gamma \mathbb{E}\mathfrak{N}^{\nu}(t), & t>0, \: \gamma > 0, \nu \in (0,1),\\
				\mathbb{E}\mathfrak{N}^{\nu}(0)=1,
			\end{cases}
		\end{align}
		where $\gamma$ is the birth intensity. The solution to \eqref{birth} reads:
		\begin{align}
			\mathbb{E}\mathfrak{N}^{\nu}(t)= E_{\nu,1}(\gamma t^{\nu}), \qquad t>0.
		\end{align}
		Recalling the definition of Mittag--Leffler function, we remark that asymptotically,
		it has a power law decay (see for example
		\citet{ogata}).
		
		Finally, a subordination relation is verified for the fractional linear death and fractional linear birth
		processes. For example, for the fractional death, we have that
		\begin{align}
			M^\nu(t) \overset{\text{d}}{=} M(V_t^\nu), \qquad t>0, \: \nu \in (0,1),
		\end{align}
		where the equality is intended for the one-dimensional distribution and $V_t^\nu$, $t>0$, $\nu \in (0,1)$, is the
		right-inverse process of the $\nu$-stable subordinator $S^\nu_x$, $x>0$, $\nu \in (0,1)$, i.e.
		\begin{align}
			V_t^\nu = \inf \{x>0 \colon S^\nu_x > t\}.
		\end{align}
	
	\section{Analytical formulation of the ETAS model}
	
		\label{analysis}
		In this section we follow \citet{sornette} to define the master self consistency equation of the ETAS model.
		We do not discuss in detail the modellistic aspects of the ETAS model (see for example \citet{helm}).
		Briefly, it can be considered as a simple branching model in which each afteshock may trigger other
		succeding aftershocks in a cascade process. It is therefore natural to regard this as an epidemic process
		in which a given \emph{parent} event of a certain magnitude, occurring at time $t=0$, gives birth
		to other \emph{child} events with rate:
		\begin{align}
			\label{1}
			\phi(t) \mathrm dt = (1-k) \theta t_0^\theta
			\frac{1}{t^{1+\theta}}H(t-t_0), \qquad t>0,
		\end{align}
		where $H(x)$, $x \in \mathbb{R}$, is the Heaviside step function,
		and $(1-k)t_0^\theta/t^{1+\theta}$ is the modified local Omori law
		on the time occurrence of the aftershocks.
		In the following, we focus our attention only on time behaviour and decay of this cascade process.
		The self consistency equation, describing the rate of seismicity at a given time $t$, reads \citep{sornette}:
		\begin{align}
			\label{ppp}
			N(t)=(1-k)\int_0^{t-t_0} N(\tau)\theta t_0 ^{\theta}(t-\tau)^{-(1+\theta)} \, \mathrm d\tau,
		\end{align}
		where $\int_{0}^{+\infty}\phi(t)\mathrm dt = (1-k)$,	
		is the branching ratio, i.e., the average number of aftershocks generated by each event;
		$t_0$ is a time delay with respect to the mainshock that occurs at time zero. Note that equation \eqref{ppp}
		is the classical Wiener--Hopf
		integral equation with a power law kernel. In \citet{sornette}, and \citet{helm}, an
		interesting discussion of the analytic solution
		to this equation was given and, in particular, three different regimes for the seismic rate $N(t)$, were found:
		\begin{itemize}
			\item $k > 0$ and $\theta > 0$: subcritical regime, less than one child per parent;
			\item $k < 0$ and $\theta>0$: supercritical regime. A transition from an Omori decay with exponent
				$p = 1-\theta$ to an explosive increase of the seismicity rate;
			\item $\theta<0$ and $k >0$: a transition from an Omori law with exponent $1-|\theta|$ to an exponentially
				increasing seismicity rate for large values of time.
		\end{itemize}
		In \citet{sornette}, a first discussion about the analytic solution of the self consistency equation was given.
		The authors found a characteristic time $t^*$,
		function of all the ETAS parameters, acting as a cut-off between the small and the large time 
		behaviour. They discussed this important feature both for the subcritical and supercitical regimes.
		For the subcritical
		regime, this means the transition 
		from an Omori exponent $p=1-\theta$ for $t<t^*$ to $p = 1+\theta$ for $t>t^*$, i.e.\
		it implies that the seismic rate
		decays slowly for small values of $t$. 
		For the supercritical regime the authors
		found a transition from an Omori power law decay with exponent $p = 1-\theta$ for
		$t<t^*$ to an explosive exponential increase. 
		In a successive paper, \citet{helm} generalised their analysis, by taking into consideration also
		the Gutenberg--Richter distribution. 
		What is more interesting for the sake of our discussion,
		is that they also analysed the third regime, i.e.\ $\theta <0$.
		This case needs more attention, because the 
		integral $ \int_{0}^{\infty}\phi(t)\mathrm dt \sim  \int_{0}^{\infty}1/(t+1)^{1+\theta} \mathrm dt$
		becomes unbounded, thus implying an infinite branching ratio.
		The authors explained this, apparently meaningless result saying that the number of children created beyond
		any time $t$ exceeds the numbers 
		of children created until time $t$.
		This interpretation is also confirmed by the discussion on the analytic solution related to this third regime,
		in which another
		charateristic time $\tau$ of transition from an anomalous slow decay for $t<\tau$
		and an explosive exponential increase for $t>\tau$, 
		similar to the supercritical case, was found.
		Therefore, this case, and the supercritical case,
		have a modellistic utility proven to be useful to understand phenomena 
		with anomalous slow seismic decay for small values of time,
		although they appear meaningless for large values of time. 
		We refer to \citet{helm} for a thorough discussion on the meaning of the regimes of the ETAS model.

		Before going ahead, we give some mathematical remarks on equation \eqref{1}.
		For example, in the last case ($\theta<0$),	we can rewrite it in the following manner:
		\begin{align}
			\label{3}
			N(t)=-(1-k)|\theta|\int_0^t N(\tau) t_0 ^{\theta}(t-\tau)^{-(1-|\theta|)} \, \mathrm d\tau, \qquad t>0,
		\end{align} 
		where we have done the approximation $t-t_0\sim t$, the time delay being small with respect to the
		total time interval.
		By recalling the definition of Riemann--Liouville fractional integral of order $\alpha>0$
		\citep{podlubny}
		\begin{align}
			J^{\alpha}f(t)=\frac{1}{\Gamma(\alpha)}\int_0^t f(\tau)(t-\tau)^{-(1-\alpha)} \, \mathrm d \tau,
			\qquad \alpha > 0, \: t>0,
		\end{align}
		with $\Gamma(\alpha)$ the Euler Gamma function, it is easy to rewrite the equation \eqref{3} as a fractional
		integral equation. Indeed we have that
		\begin{equation}
			\label{prefrac}
			N(t)=-\lambda J^{|\theta|} N(t), \qquad t>0, \: \theta<0,
		\end{equation}
		where $\lambda= (1-k)|\theta|\Gamma(|\theta|)$.
		Finally, we arrive at the following fractional differential equation:
		\begin{align}
			\label{fracc}
			D^{|\theta|}N(t)=-\lambda N(t), \qquad t>0, \: \theta < 0,
		\end{align}
		where $D^{|\theta|}$ is the Caputo fractional derivative of order $|\theta| > 0$, defined as
		\citep{podlubny}
		\begin{align}
			\label{caputo}
			D^{|\theta|} N(t) = \frac{1}{\Gamma(m-|\theta|)} \int_0^t \frac{\frac{\mathrm d^m}{
			\mathrm d \tau^m}N(\tau)}{(t-\tau)^{|\theta|-m+1}}
			\mathrm d\tau,
		\end{align}
		where $m=\lceil |\theta| \rceil$.
		
		Equation \eqref{prefrac} can thus be written as a
		simple linear fractional differential equation. Equation \eqref{fracc},
		in the case $\theta \in (-1,0)$, with the initial condition $N(0)=n_0$,
		$n_0 \in \mathbb{N}$ (meaning that the initial
		rate of seismicity is $n_0$), can be easily solved as we have shown in the previous section.
		
		Finally, we realize that the ETAS self consistency equation
		is formally similar to equation \eqref{ponte} governing the mean behaviour
		of the fractional linear death process. Furthermore, the supercritical regime with explosive growth
		can be linked in the same way to the fractional linear pure birth process.
		In the next section, we discuss similarities and differences 
		between these models.	

	\section{Relations between self consistency equation of the ETAS model and fractional birth or death processes} 

		\label{sec4}
		Taking inspiration from the formal analogy shown in Section \ref{analysis},
		here, we discuss the analytical behaviour of
		the rate $N(t)$, $t>0$, in the framework of fractional death and fractional birth processes,
		considering, on the one hand, the behaviour predicted by them,
		and, on the other hand, the physical meaning of these predictions in the framework of the ETAS model.
		In short, the rate of seismicity $N(t)$, $t>0$, is considered here as the mean behaviour of an
		underlying fractional linear death or birth process.
	
		Similarly to the ETAS model, starting from equation \eqref{fracc}
		describing the rate of an epidemic process, we can distinguish two regimes:
		\begin{itemize}
			\item if $\lambda>0$, we have
				\begin{align}
					\label{fraccc}
					D^{\nu}N(t)=-\lambda N(t), \qquad t>0,\;\nu \in (0,1).
				\end{align}
				As we have shown, it is identical to
				the fractional differential equation arising in the fractional death process.  
				Physically it corresponds to the subcritical regime in the framework of the ETAS model.
				However we must underline some differences, remembering that our formal analogy is not rigorously
				an identity. We observe
				a strong exponential decay for small values of $t$ and an asymptotic power law decay
				when $t \rightarrow \infty$.
				Thus, we have a reasonable physical picture for cases
				with a slow decay for large times, but this prediction is
				slightly different from that of ETAS.
			
				In Fig.\ \ref{mea} the solution to the Cauchy problem composed by \eqref{fraccc}
				and the initial condition $N(0) = n_0 = 1$, for some values of $\nu$,
				is shown; it corresponds to the mean value of
				the linear fractional death process. Note the fast decay for small values of $t$.
				\begin{figure}
					\centering
					\includegraphics[scale=0.8]{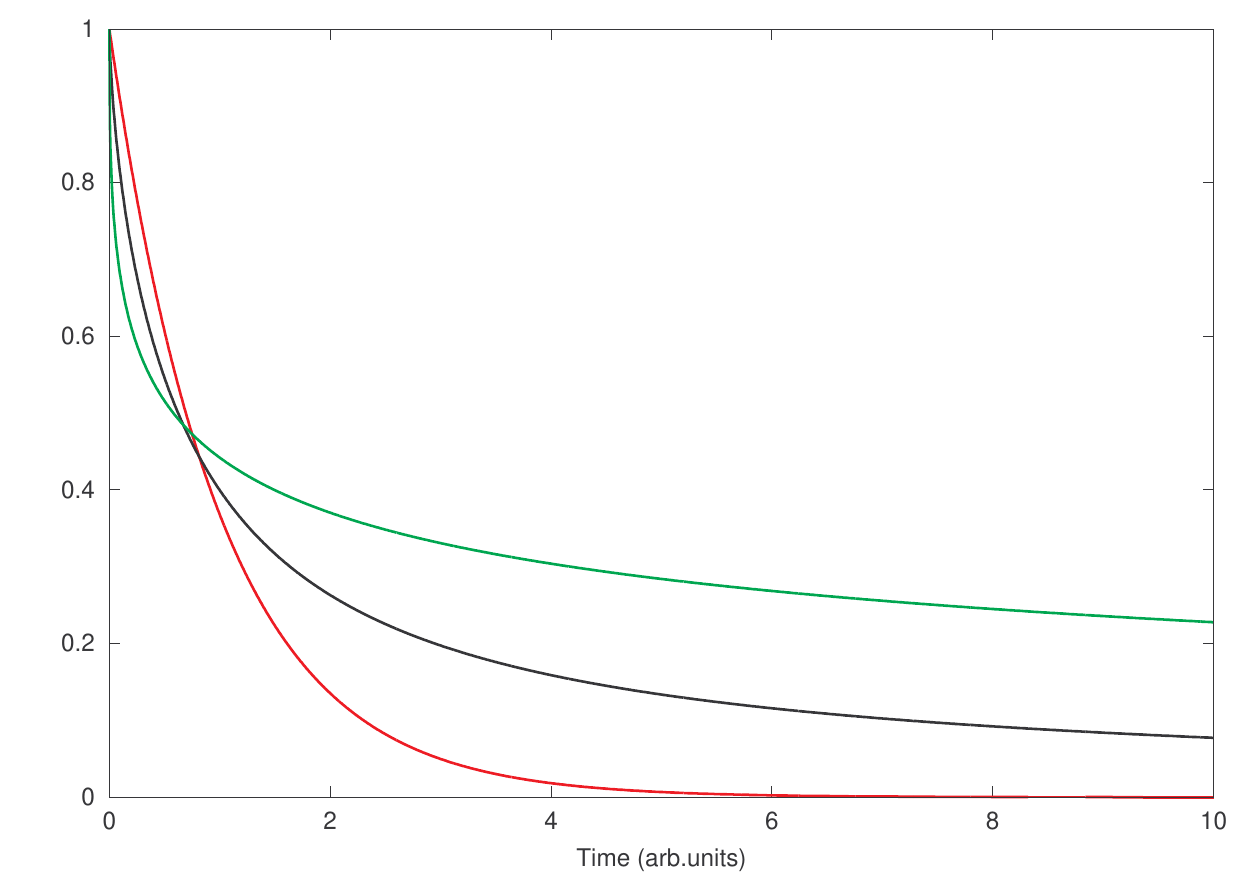}
					\caption{The mean value of the linear fractional death process with $\nu=1$ (red), $\nu=0.7$ (black)
					and $\nu=0.4$ (green). \label{mea}}				
				\end{figure}

				\item if $\lambda<0$ we have the following equation:
				\begin{align}
					\label{wicklow}
					D^{\nu}N(t)=|\lambda| N(t), \;\nu \in (0,1).
				\end{align}
				This corresponds to the fractional differential equation arising in the linear fractional birth process.
				This case is equivalent to the supercritical regime of the ETAS model.
				Indeed, we have an explosive growth. Fig.\ \ref{meanpol}
				presents the solution to the Cauchy problem composed by \eqref{wicklow} and the
				initial condition $N(0)=1$. We have an exponential growth,
				as for the supercritical regime of the ETAS model.
				\begin{figure}
					\centering
					\includegraphics[scale=0.8]{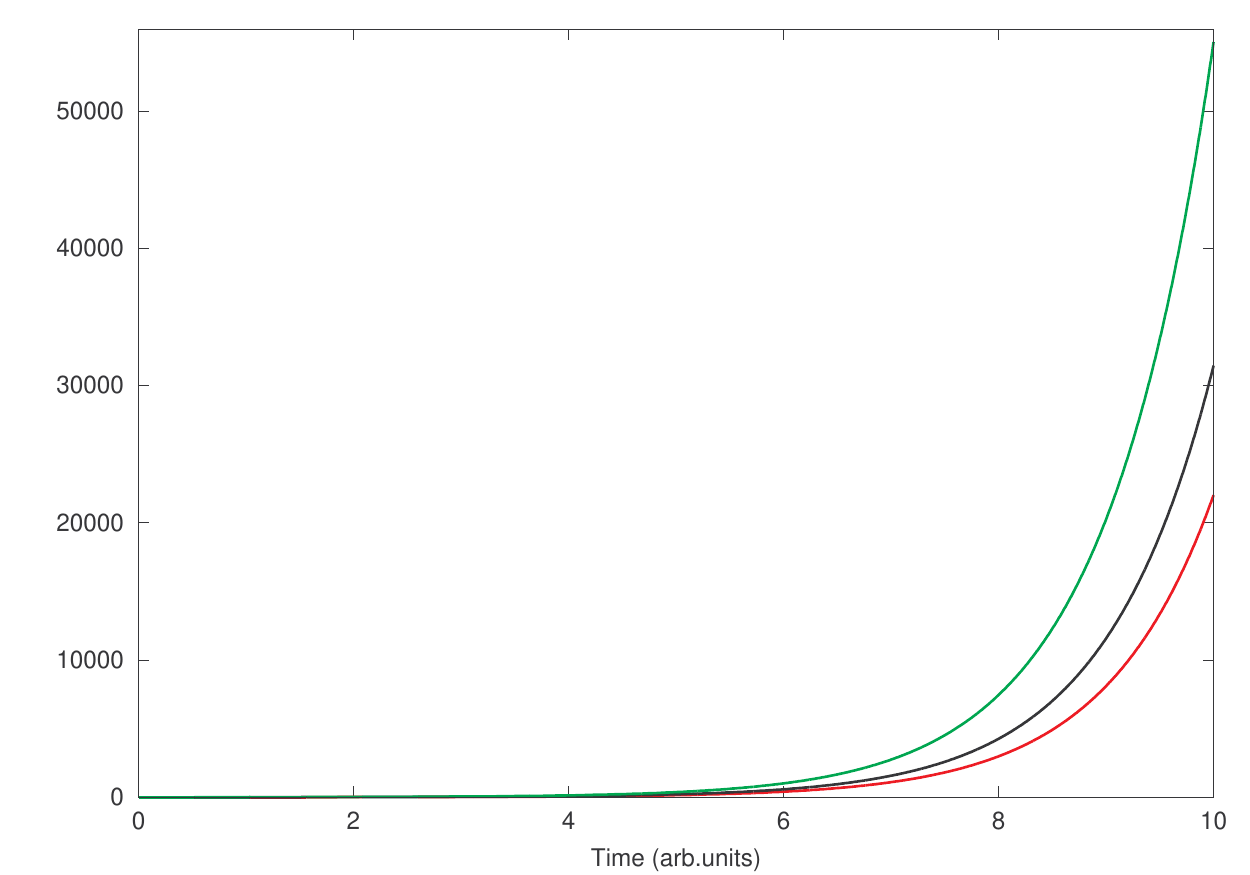}
					\caption{The mean value of the linear fractional birth process with $\nu=1$ (red), $\nu=0.7$ (black)
					and $\nu=0.4$ (green). \label{meanpol}}				
				\end{figure}
		\end{itemize}

	\section{An application of the fractional pure birth process: theoretical model of tumoral growth}

		\label{sec5}	
		Beyond the formal analogy with the ETAS model we can explain the utility of our proposal
		with a simple toy model about tumoral growth.
		In the following we show that the fractional linear pure birth process
		can give an interesting framework to modelling this kind of phenomena and have good
		agreement with some results in this field.
		Assume that we have an initial number of infected cells $\mathfrak{N}^\nu(0)=n_0$.
		We are interested in building a theoretical model about the growth
		of cancer cells with respect to time in the dynamics of metastasis.
		Suppose that this process is described by a linear fractional pure birth model. 
		Therefore we have (see \citet{pol})
		that the distribution of the population size of cancer cells at time $t$ is given by:
		\begin{align}
			\label{n0}
			& \text{Pr} \{ \mathfrak{N}^\nu(t) = k+n_0 | \mathfrak{N}^\nu(0) = n_0 \} \\
			& = \binom{n_0+k-1}{k} \sum_{r=0}^{k}
			\binom{k}{r}(-1)^r E_{\nu,1}(-(n_0+r)\gamma t^{\nu}), \quad t>0, \: \nu \in (0,1], \: \gamma>0. \notag
		\end{align}
		Now, from \eqref{n0}, we can infer the probability of a new offspring at the beginning of the process:
		\begin{align} 
			& \text{Pr} \{ \mathfrak{N}^\nu(\mathrm dt)=1+n_0|\mathfrak{N}^\nu(0)=n_0 \} \\
			& = n_0 \left[ E_{\nu,1}(-n_0\gamma (\mathrm dt)^{\nu})-E_{\nu,1}(-\gamma(n_0+1)
			(\mathrm dt)^{\nu}) \right] \sim n_0 \frac{\gamma(\mathrm dt)^{\nu}}{\Gamma(\nu+1)}, \notag
		\end{align}
		by writing only the lower order terms. This result shows that the probability of a new offspring is proportional
		to the time interval $(\mathrm dt)^{\nu}$ and to the initial number of progenitors.
		This is an interesting picture that can be realistic
		for a large number of complex growth dynamics with emerging power law behaviour. 
		From the previous analysis, we also have that the average expansion of the population
		of cancer cells is given by:
		\begin{align}
			\mathbb{E} \mathfrak{N}^\nu(t)= n_0 E_{\nu,1}(\gamma t^{\nu}), \qquad t>0.
		\end{align}
		In conclusion we have a theoretical framework that can be used to describe growth processes
		exhibiting power law expansions.
		
		We discussed this model by referring to the cancer
		growth dynamics because several works (see for example \citet{dattoli,west})
		studied a similar dynamic starting from the Kleiber law.
		It is not the purpose of this paper going inside this vaste
		field of study, but we remark that we have a conceptual well posed model without
		empirical assumptions. Thus we obtain a probabilistic 
		point of view that is more \textit{fundamental} with respect to the dynamic of these processes. 
		Finally, we must observe that, in order to make the model more realistic,
		it is possible to directly generalise the fractional birth process by introducing a
		saturation threshold. This will be matter of a future paper.
		
		We also notice that the role played by the allometric coefficients
		in the Kleiber law, in our model, is played by the real order of derivation $\nu$.
		These allometric coefficients are empirically established;
		we have the same picture but with only one free fitting parameter.
		Beyond the probabilistic clear view of the processes, this is also a strong point of this theory,    

	\section{Discussion}
		
		\label{sec6}
		We now move to understand the differencess between the ETAS model and the fractional epidemic processes
		predictions. First, we must identify the real order of derivation $\nu$ with the parameter
		$\theta$ of the Omori law 
		and $\lambda$ with $\lambda= (1-k)|\theta|\Gamma(|\theta|)$ of the ETAS model. Although we do not have
		exactly the same cut-off 
		between different regimes as in the ETAS model, we obtain some reasonable agreement
		for various ranges of the parameters.
		Moreover, in \citet{helm}, by discussing the self-consistency equation for the supercritical regime,
		a cut-off between slow decay for small times and explosive growth was found.
		In our case, the supercritical regime is less meaningful, presenting immediately an exponential growth.
		On the other hand, from the self-consistency equation \eqref{ppp}
		and by considering the physical behaviour in the subcritical regime,
		the agreement of the two models is clear, both describing epidemic processes.
	
		At this point, by keeping in mind the relations between the ETAS model and the fractional linear
		birth or death processes, we make some remarks.
		We have a generalised analytic model with clear relations with epidemic models.
		We also notice that in many epidemic processes in complex systems (as in seismology and biology),
		empirical power law behaviour emerges. We have seen that the prediction on the rate of fractional processes,
		i.e., the Mittag--Leffler function (generalised exponential), takes into account this behaviour asymptotically.
  		We conclude that the fractional processes we treated in this paper, are an important instrument
  		in order to study
  		epidemic processes which exhibit
  		a rapid decay for small times and slow asymptotic power law decay for large times.

	\bibliographystyle{abbrvnat}
	\bibliography{ETAS4}
	\nocite{*}

\end{document}